\documentclass[12pt,a4paper]{article}
\usepackage[latin2]{inputenc}
\usepackage[mathscr]{eucal}
\usepackage{t1enc}
\usepackage{amssymb}
\usepackage{amsfonts}
\usepackage{amsmath}

\usepackage{amsthm}
\usepackage{latexsym}
\usepackage[dvips]{color}

\usepackage[dvips]{epsfig}

\usepackage{bbm}

\newtheorem{Theorem}{Theorem}[section]
\newtheorem{Remark}[Theorem]{Remark}
\newtheorem{Lemma}[Theorem]{Lemma}

\newtheorem{Proposition}[Theorem]{Proposition}

\newtheorem{Corollary}[Theorem]{Corollary}

\begin{document}

\title{
{\bf The high temperature Ising model \\ on the triangular lattice \\ is a critical percolation model}}

\author{
{\bf Andr\'as B\'alint},
\thanks{E-mail: abalint@few.vu.nl}
{\bf \, Federico Camia},
\thanks{Research supported in part by a Veni grant of the NWO
(Dutch Organization for Scientific Research).}\,
\thanks{E-mail: fede@few.vu.nl}
{\bf \, Ronald Meester}
\thanks{Research supported in part by a Vici grant of the NWO
(Dutch Organization for Scientific Research).}\,
\thanks{E-mail: rmeester@few.vu.nl}\\
{\sl Department of Mathematics, Vrije Universiteit Amsterdam}
}

\maketitle

\begin{abstract}
The Ising model at inverse temperature $\beta$ and zero external field 
can be obtained via the Fortuin-Kasteleyn (FK) random-cluster model with $q=2$ and density
of open edges $p=1-e^{-\beta}$ by assigning spin $+1$ or $-1$ to each vertex in such a way
that (1) all the vertices in the same FK cluster get the same spin and (2) $+1$ and $-1$
have equal probability.

We generalize the above procedure by assigning spin $+1$ with probability
$r$ and $-1$ with probability $1-r$, with $r \in [0,1]$, while keeping condition (1).
For fixed $\beta$, this generates a dependent (spin) percolation model with
parameter $r$. We show that, on the triangular lattice and for $\beta<\beta_c$,
this model has a percolation phase transition at $r=1/2$, corresponding
to the Ising model. This sheds some light on the conjecture that the high
temperature Ising model on the triangular lattice is in the percolation
universality class and that its scaling limit can be described in terms of SLE$_6$.

We also prove uniqueness of the infinite $+1$ cluster for $r>1/2$, sharpness of
the percolation phase transition (by showing exponential decay of the cluster
size distribution for $r<1/2$), and continuity of the percolation function for
all $r \in [0,1]$.
\end{abstract}


\section{Introduction and motivation}

If one considers the percolation properties of spin clusters, the high temperature
($\beta<\beta_c$) Ising model on the {\em triangular lattice} $\mathbb T$ with no
external field is believed to show critical behaviour and to be in the universality
class of Bernoulli (independent) percolation. One way to understand the conjecture
is in terms of the renormalization group: one expects the high temperature phase of
the Ising model to be in the basin of attraction of the (stable) infinite temperature
fixed point, which in the case of the triangular lattice (but {\em not} of other
lattices) corresponds to critical site percolation.

The conjecture can also be cast in the language of the $O(n)$ model. Indeed, the high
temperature phase of the zero field Ising model on $\mathbb T$ falls in the so-called
dense phase of the $O(n)$ model with $n=1$, which is believed to show critical behaviour
and to be driven, under the action of the renormalization group, to a stable fixed point
corresponding, once again, to critical percolation.
Rephrased in the language of scaling limits, the conjecture can be expressed in terms of
convergence of certain interfaces to SLE$_6$, and appears in several places (see, e.g.,
\cite{KagerNienhuis,Sheffield,smirnov,BCM}),
often as part of a larger
conjecture about the $O(n)$ model.

Motivated by the above considerations, in this paper we show that, for any fixed inverse
temperature $\beta<\beta_c$, the zero field Ising model on $\mathbb T$ corresponds to
the critical point of a dependent percolation model. To do that, we generalize the FK
random-cluster representation of the Ising model, obtaining a site percolation model
with parameter $r \in [0,1]$ which reduces to the Ising model for $r=1/2$. We then show
that the new percolation model has a (sharp) percolation phase transition at $r=1/2$.

It will be clear from its definition that in the new percolation model the states
of different sites are not independent, but that correlations decay exponentially
with the distance. As a consequence, we can now view the conjecture about the high
temperature Ising model on the triangular lattice as a particular instance of another
general conjecture that has to do with percolation only, namely that the scaling
limit of a (two-dimensional) percolation model exists and is independent of the
particular model, as long as the correlations in the probability measure decay fast
enough (see, e.g., \cite{Cardy}).
Such a result has been proved for a few specific models of dependent percolation
(see \cite{Camia-Newman-Sidoravicius,camia-newman,Camia-bootstrap,Camia-enhancement,Binder-Chayes-Lei}).

It is interesting to notice that the Ising model corresponds to the self-dual point
of the new percolation model. Therefore, our result also shows that for that model
the self-dual and the critical point coincide, in accordance with a very natural
principle which is believed to be valid in great generality, but which has been
verified only in a handful of cases, including bond percolation on the square
lattice~\cite{Kesten}, site percolation (see~\cite{kesten-book}) and the Divide
and Colour (DaC) model~\cite{BCM} on the triangular lattice, and Voronoi
percolation~\cite{Voronoi}. The same priciple should apply to other interesting models,
such as the random-cluster model (where it is known for $q=1$, corresponding to
percolation, and $q=2$, corresponding to the Ising model --- see~\cite{grimmett2}),
other DaC models (see~\cite{BCM}, and in particular Conjecture~1.7 there) and
``confetti percolation'' (see Problem 5 in~\cite{bs}).

In our analysis we make substantial use of a result by Higuchi, who has extensively
studied the percolation properties of the two-dimensional Ising model 
(see \cite{HiguchiRSW,Higuchiperc,Higuchiremark,Higcoex,Higuchi}).

\section{Main results}\label{section-main-results}

We work on the triangular lattice $\mathbb{T}$ with vertex set $\mathcal{V}_{\mathbb{T}}$
and edge set $\mathcal{E}_{\mathbb{T}}$, and denote the unique {\em Ising Gibbs measure} on
$\mathbb{T}$ at inverse temperature $\beta < \beta_c$ and zero external field by $\mu_{\beta}$.

The {\em random-cluster measure} $\nu_{p,q}$ on edge configurations
$\eta \in \{0,1\}^{\mathcal{E}_{\mathbb{T}}}$  (with the usual $\sigma$-field generated by
cylinder events) is characterized by two parameters satisfying $0 \leq p \leq 1$ and $q>0$
(see~\cite{grimmett2} for the definition and some background).
We call an edge $e\in \mathcal{E}_{\mathbb{T}}$ \textit{open} if $\eta(e)=1$, and
\textit{closed} otherwise. The maximal connected components of the graph obtained
by removing all the closed edges from $\mathbb{T}$ are called \textit{FK clusters}.

For fixed $q$, the random-cluster measure has a percolation phase transition at some
$0<p_c(q)<1$, and with probability one all FK clusters are finite if $p < p_c(q)$.
When $q=2$, one can generate an Ising spin configuration
$\sigma \in \{+1,-1\}^{\mathcal{V}_{\mathbb{T}}}$ distributed according to $\mu_{\beta}$,
$\beta<\beta_c$, by drawing an edge configuration according to $\nu_{p,2}$ with
$p=1-e^{-\beta}$ and assigning spin $+1$ or $-1$ to each vertex of $\mathbb T$ in
such a way that (1) all the vertices in the same FK cluster get the same spin and
(2) $+1$ and $-1$ have equal probability. Note that $\beta<\beta_c$ implies
$p=1-e^{-\beta}<p_c(2)$, so that the FK clusters are all finite with probability one.

From now on, unless otherwise stated, we will always assume that $q=2$.
We generalize the above procedure by assigning spin $+1$ with probability $r$
and $-1$ with probability $1-r$, with $r \in [0,1]$, while keeping condition (1).
For fixed $\beta$, this generates a dependent (spin) percolation model with
parameter $r$, whose measure we denote by ${\mathbb P}_{\beta,r}$. Clearly,
the spin marginal of ${\mathbb P}_{\beta,1/2}$ coincides with $\mu_{\beta}$.
Note also that ${\mathbb P}_{0,1/2}$ (equivalently, $\mu_0$) is a product measure
and corresponds to critical site percolation on $\mathbb T$. As soon as $\beta>0$,
however, the spins are correlated. Nonetheless, the exponential decay of the FK
cluster size distribution when $\beta<\beta_c$ (see~\cite{grimmett2}) immediately
implies the exponential decay of correlations in the measure ${\mathbb P}_{\beta,r}$.

We call a maximal connected subset $V$ of ${\mathcal{V}_{\mathbb{T}}}$ such
that all vertices in $V$ have the same spin a {\em spin cluster}. If the spins
in $V$ are all $+1$ (respectively, $-1$), we call $V$ a {\em ($+$)-cluster}
(resp., a {\em ($-$)-cluster}). Our aim is to study the percolation properties
of spin clusters. We denote by $\Theta(\beta,r)$ the $\mathbb{P}_{\beta,r}$-probability
that a given vertex of the triangular lattice is contained in an infinite ($+$)-cluster,
and define $r_c(\beta):=\sup \{r:\Theta(\beta,r)=0\}$. By the size of a cluster
we mean the number of vertices in the cluster. The main result of this paper
is the following theorem where, due to the $+/-$ symmetry of the model, we
focus without loss of generality on the behaviour of ($+$)-clusters.
\begin{Theorem} \label{rc12}
For all $\beta<\beta_c$, $r_c(\beta)=1/2$. Moreover,
\begin{itemize}
\item If $r<1/2$, the distribution of the size of the ($+$)-cluster of the origin
has an exponentially decaying tail.
\item If $r=1/2$, $\Theta(\beta,1/2)=0$ and the mean size of the ($+$)-cluster of
the origin is infinite.
\item If $r>1/2$, there exists a.s.\ a unique infinite ($+$)-cluster.
\end{itemize}
\end{Theorem}

Note that $r=1/2$ is clearly the self-dual point of ${\mathbb P}_{\beta,r}$.
Thus, Theorem~\ref{rc12} implies that the critical point of the model
coincides with its self-dual point. We remark that  
one can obtain a polynomial lower bound for the 
tail distribution of the ($+$)-cluster of the origin at $r=1/2$
by using elementary duality arguments only, see \cite{GrimmettJanson}, p.\ 15.

Since the phase transition described in Theorem~\ref{rc12} is continuous,
one may expect continuity of the percolation function $\Theta(\beta,r)$.
Indeed, this can be proved by standard methods.

\begin{Theorem} \label{contbond}
For each $\beta<\beta_c$, $\Theta (\beta,r)$ is a continuous function of
$r \in [0,1]$.
\end{Theorem}

It is worth remarking that ${\mathbb P}_{\beta,r}$ belongs to a family of
measures that can be obtained via a two step procedure: first partition the
vertices of a lattice into clusters according to some rule, then assign
spin values (or colours) to the vertices with some probability, making sure
that all vertices in the same cluster get the same colour. We call such
measures {\em Divide and Colour (DaC) models}. The first DaC model was
introduced by H\"aggstr\"om~\cite{DaC}, and its phase transition is studied
in detail in~\cite{BCM}. DaC models can be considered as natural dependent
percolation models. They are relatively simple, yet their analysis is
considerably more complicated than that of Bernoulli percolation (see,
e.g., the present paper and~\cite{BCM}), and requires new
techniques that may be useful in studying other dependent models.

A brief outline of the paper is given as follows.
In Section \ref{preliminaries}, we introduce some more definitions and notation,
and we collect results which are either known or can be proved by standard methods,
including a result by Higuchi~\cite{Higcoex} about the Ising model.
We shall use them later, together with the standard Edwards-Sokal coupling \cite{ES} and
results described in Section \ref{domination-lemmas}, which contains some technical
lemmas and an overview of the main step in the proof of Theorem \ref{rc12}. In
Section~\ref{bondmainproofsec}, we prove Theorems \ref{rc12} and \ref{contbond}.

\section{Preliminaries}\label{preliminaries}

\subsection{Notation and definitions}\label{defs}

\begin{figure}[h]
\centering
\includegraphics[scale=0.37, trim= 0mm 0mm 0mm 0mm, clip ]{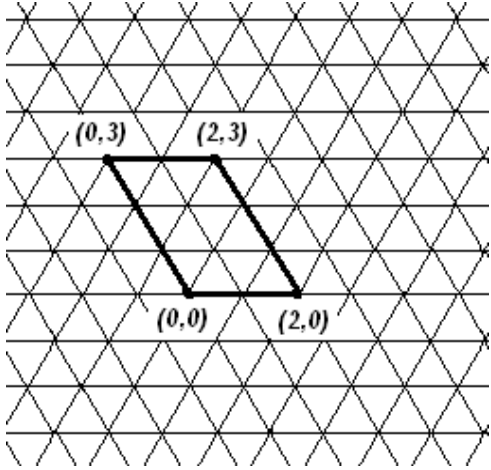}
\caption{Portion of the triangular lattice $\mathbb{T}$. The heavy segments are the sides of
the parallelogram $S_{2,3}=[0,2]\times [0,3]$.}
\label{rec}
\end{figure}

In order to define a concrete coordinate system in the triangular lattice $\mathbb T$,
we embed $\mathbb T$ in ${\mathbb R}^2$ as in Figure~\ref{rec}, so that its set of vertices $\mathcal{V}_{\mathbb T}$
consists of the intersections of the lines $y=-\hspace{0.03cm}\sqrt{3}\hspace{0.05cm} x+\sqrt{3}\hspace{0.05cm} k$
and $y=\frac{\sqrt{3}}{2}\hspace{0.05cm} \ell$ for $k,\ell \in \mathbb{Z}$, and denote the elements of
$\mathcal{V}_{\mathbb T}$ by $(k,\ell)$.
We call two vertices in $\mathcal{V}_{\mathbb T}$ \textit{adjacent} if their Euclidean distance is 1,
and define the edge set $\mathcal{E}_{\mathbb{T}}$ by $\mathcal{E}_{\mathbb{T}}=\{(v,w):v$ and $w$ are adjacent$\}$.

The state space of our configurations is denoted by $\Omega :=\Omega _{D}\times \Omega _{C}$,
where $\Omega _{D}=\{0,1\}^{\mathcal{E}_{\mathbb T} }$ is the set of random-cluster realisations,
and $\Omega _{C}=\{-1,+1\}^{\mathcal{V}_{\mathbb T}}$ corresponds to the spin configurations.
The probability measure $\mathbb{P}_{\beta ,r}$ is the measure (on the usual $\sigma$-algebra on $\Omega$)
obtained by the procedure described in Section \ref{section-main-results};
we denote the expectation with respect to $\mathbb{P}_{\beta ,r}$
with $\mathbb{E}_{\beta ,r}$.

We introduce the set $\tilde\Omega \subset \Omega$ as the set of configurations such that
vertices in the same FK cluster have the same spin, and we equip $\tilde\Omega $ with a
partial order, which depends on the spins only, as follows.
For $\omega _1=(\eta _1,\sigma _1), \omega _2=(\eta _2,\sigma _2) \in \tilde\Omega $ we say that
$\omega _1 \geq \omega _2$ if $\sigma _1 (x)\geq \sigma _2(x)$ holds for every $x\in \mathbb{T}$.
All the configurations are implicitly assumed to be in $\tilde{\Omega}$.
We call an event $A \subset \tilde{\Omega}$ \textit{increasing} if $\omega \in A$ and
$\omega '\geq \omega $ implies $\omega '\in A$.
$A$ is a \textit{decreasing event} if $A^c$ is increasing.

We call a sequence $(x_0,x_1,\ldots ,x_n)$ of vertices in $\mathbb{T}$ a \textit{(self-avoiding) path}
if for all $i=0,\ldots ,n-1$, $x_i$ and $x_{i+1}$ are adjacent, and for any
$0\leq i< j\leq n,\hspace{0.1cm} x_i\neq x_j$.
A \textit{horizontal crossing} of a parallelogram $R=[a,b]\times [c,d]$, with $a,b,c,d \in {\mathbb Z}$,
is a path ${x_0,x_1,\ldots ,x_n}$ such that $x_0\in \{a\}\times [c,d]$, $x_n\in \{b\}\times [c,d]$
and for all $i$, $x_i\in R$.
A \textit{vertical crossing} of the same parallelogram is a path ${x_0,x_1,\ldots ,x_n}$ such
that $x_0\in [a,b]\times \{ d\}$, $x_n\in [a,b]\times \{c\}$ and for all $i$, $x_i\in R$.

In a configuration $(\eta ,\sigma )\in \tilde \Omega $, a \textit{($+$)-path} is a path ${x_0,x_1,\ldots ,x_n}$
such that for all $i=0,\ldots ,n, $ $\sigma(x_i)=+1$.
\textit{Horizontal ($+$)-crossings} and \textit{vertical ($+$)-crossings} are defined
analogously.
The definitions of \textit{($-$)-path}, \textit{horizontal ($-$)-crossing},
\textit{vertical ($-$)-crossing} are obtained by replacing $+1$ with $-1$.

Let $S_{n,m}$ denote the parallelogram $[0,n]\times [0,m]$, with $n,m \in {\mathbb N}$.
Denote by $V^+_{n,m}$ the event that there is a vertical ($+$)-crossing in $S_{n,m}$;
let $H_{n,m}^+$ be the corresponding event with a horizontal ($+$)-crossing.
The analogous events with ($-$)-crossings are denoted by $V_{n,m}^-$ and $H_{n,m}^-$,
respectively.

Let $d$ denote the graph distance on $\mathbb{T}$. We define the distance between two sets $V$ and $W$ by
$d(V,W)=\{\min (d(v,w)):v\in V,w\in W\}$.
Let $B(v,n)$ denote the disc of radius $n$ with center at vertex $v$ in the metric $d$,
i.e., $B(v,n)=\{w:d(v,w)\leq n\}$.
For a vertex set $A\subset \mathcal{V}_{\mathbb{T}}$, we denote by $\partial A$ the \textit{vertex boundary} of $A$,
that is, we define $\partial A=\{v\in A:\exists \hspace{0.05cm} w\in {\cal V}_{\mathbb{T}}\setminus A $ such that
$d(v,w)=1 \}$.
For a vertex $v\in \mathcal{V}_{\mathbb{T}}$, let $C_v ^{FK}$ be the FK cluster of $v$, i.e., the set of
vertices that can be reached from $v$ through edges that are open in the underlying random-cluster measure
with parameters $p$ and $q=2$.
Let us define the \textit{dependence range} of a vertex $v$ by
$\mathcal{D}(v):=\max \{n\in \mathbb{N}:C_v ^{FK}\cap \partial B(v,n) \neq \emptyset\}$.

We call an edge set $E=\{e_1,e_2,\ldots ,e_k\}$ a \textit{barrier} if
removing $e_1,e_2,\ldots ,e_k$ (but not their end-vertices) separates the graph $\mathbb{T}$ into two
or more disjoint connected subgraphs, of which exactly one is infinite.
(Note that a barrier as defined above corresponds to dual circuits in bond percolation. Its definition
is motivated by Lemma \ref{conditional-independence}.)
We call the infinite subgraph 
the \textit{exterior} of $E$, and denote it by $ext(E)$.
We call the union of the finite subgraphs the \textit{interior} of $E$, and denote it by $int(E)$.
With an abuse of notation, we shall write $int(E)$ and $ext(E)$ also for the vertex sets
of $int(E)$ and $ext(E)$ whenever it does not cause confusion.
$E=\{e_1,e_2,\ldots ,e_k\}$ is a \textit{closed barrier} in a configuration $(\eta ,\sigma )\in \tilde\Omega $
if $E$ is a barrier and $\eta (e_i)=0$ holds for $i=0,\ldots ,k$.
For a vertex set $A\subset \mathcal{V}_{\mathbb{T}}$, let $\Delta A$ denote the \textit{edge boundary} of $A$, that is,
$\Delta A=\{(x,y)\in \mathcal{E}_{\mathbb{T}} : x\in A, y\in \mathcal{V}_{\mathbb{T}} \setminus A\}$.
Note that for $\beta < \beta_c$, the edge boundary of any FK cluster is a.s.\ a closed barrier.

\subsection{Preliminary results}\label{bond-preliminary-results}

To make the paper self-contained, we collect here the tools needed to prove Theorems \ref{rc12}
and \ref{contbond}. The first theorem in this subsection follows from results in \cite{abf},
and is stated explicitly e.g.\ in \cite{grimmett2}.

\begin{Theorem}
\label{lexpdecay2}
If $p<p_c(2)$, there exists $\psi (p)>0$ such that for all n, we have
\begin{displaymath}
\nu _{p,2}(\mathcal{D}(0)\geq n)\leq e^{-n\psi (p)}.
\end{displaymath}
\end{Theorem}

Another property of the random-cluster measures is that
for $e\in \mathcal{E}_{\mathbb{T}}$ the conditional measure $\nu _{p,q}(\cdot \mid Y(e)=0)$
can be interpreted as a random-cluster measure with the same parameters $p$ and $q$ on
the graph obtained from $\mathbb{T}$ by deleting $e$ (see \cite{grimmett2}).
This property implies the following observation, which we state as a lemma for ease of reference.
\begin{Lemma}\label{conditional-independence}
If $B=\{e_1,\ldots ,e_k\}$ is a barrier,
$C(B)=\cap _{i=1}^k\{\eta (e_i)=0\}$,
$E_{1}$ and $E_{2}$ are events which depend only on states of edges and spins
of vertices in $int(B)\cup B$ and $ext(B)$ respectively, then conditioned on $C(B)$,
$E_1$ and $E_2$ are independent.
\end{Lemma}

In the proof of Theorem \ref{rc12}, we will use a version of Russo's formula for decreasing events,
hence we state the theorem in a slightly unusual form.
The proof, as sketched in \cite{BCM}, is standard.
Let $A$ be an event, and let $\omega =\left( \eta ,\sigma \right) $ be a configuration
in $\tilde\Omega$. Let $C$ be an FK cluster in $\eta $.
We call $C$ \textit{pivotal} for the pair $\left( A,\omega \right) $ if $I_A(\omega )\neq I_A(\omega ')$
where $I_A$ is the indicator function of $A$, $\omega' = (\eta,\sigma')$, and $\sigma'$ agrees with $\sigma$
everywhere except that the spins of the vertices in $C$ are different.

\begin{Theorem}
\label{Russocons}
Let $W$ be a set of vertices with $|W|<\infty $, and let $A$ be
a decreasing event that depends only on the spins of vertices in $W$.
Then we have that
\begin{displaymath}
\frac{d}{dr}\mathbb{P}_{\beta ,r}(A)=-\mathbb{E}_{\beta ,r}(n(A)),
\end{displaymath}
where $n(A)$ is the number of FK clusters which are pivotal for $A$.
\end{Theorem}

The following result, like Lemma 2.10 in \cite{BCM}, is a finite size criterion for percolation.

\begin{Lemma}\label{finitesizecrit2}
There exists a constant $\varepsilon >0$ with the
following property. If $\beta , p=1-e^{-\beta }$ and $N\in \mathbb{N}$ satisfy
\begin{equation}\label{littlerange}
(N+1)(3N+1)\nu _{p,2}(\mathcal{D}(0) \geq \frac{N}{3}) \leq \varepsilon \nonumber
\end{equation}
and
$$
\mathbb{P}_{\beta ,r}(V_{N,3N}^+) >  1-\varepsilon ,
$$
then $\Theta (\beta ,r)>0$.
\end{Lemma}

As in \cite{BCM}, this theorem can be proved by a
coupling argument with a 1-dependent bond percolation model.
Theorem \ref{lexpdecay2} and Lemma \ref{finitesizecrit2} imply the following result.

\begin{Theorem}\label{implperc}
For all $\beta <\beta _c$, if
$\limsup \limits _{n\to \infty } \mathbb{P}_{\beta ,r}(V_{n,3n}^+)=1$ for some $r$,
then $\Theta (\beta ,r)>0$.
\end{Theorem}

\subsection{Cut points}\label{Isingsub}
\begin{figure}[h]
\centering
\includegraphics[scale=0.34, trim= 0mm 0mm 0mm 0mm, clip ]{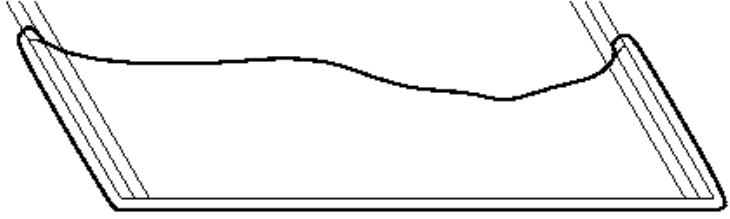}
\caption{Bottom part of the parallelogram $S_{n,6n}$.
The lines inside $S_{n,6n}$ parallel to the sides of $S_{n,6n}$
represent the sides of the parallelograms $S_{n,6n}^{\prime }$
and $S_{n,6n}^{\prime \prime }$.
The partly thin, partly thick line from the left side of $S_{n,6n}$ to the right side
represents a horizontal crossing $R\in {\cal R}_{n,4n}$.
The thick line represents the boundary of $D(R)$.}
\label{bottompar}
\end{figure}

We shall use (a slightly modified version of)
a result of Higuchi from \cite{Higcoex} (see also Proposition 4.2 in \cite{Higuchi}) about the Ising model.
In order to state the theorem, we need a few definitions. For positive integer values of $k$,
let ${\cal R}_{n,kn}$ be the collection of all horizontal crossings in $S_{n,kn}$. For $R\in {\cal R}_{n,4n}$,
we denote the region in $S_{n,6n}$ (note the different side length)
under $R$ by $L(R)$, the region in $S_{n,6n}$ above $R$ by $A(R)$, the parallelogram
$[\lfloor n^{1/4}\rfloor ,n-\lfloor n^{1/4}\rfloor ]\times [0,6n]$ by $S_{n,6n} ^{\prime }$, and the parallelogram
$[2\lfloor n^{1/4}\rfloor ,n-2\lfloor n^{1/4}\rfloor ]\times [0,6n]$ by $S_{n,6n} ^{\prime \prime }$.
(For $a\in \mathbb{R}$, we denote by $\lfloor a\rfloor $ the greatest integer smaller than or equal to $a$.)
Also, let $D(R)$ denote the vertex set
$\{v\in \mathcal{V}_{\mathbb{T}}\setminus (A(R)\cap S_{n,6n}^{\prime }): d(v,L(R)\cup R)\leq n^{1/4}\}$
(see Figure \ref{bottompar}).
We call a vertex $x\in R$ a \emph{cut point of $R$ in $S_{n,6n}$}
if there exists a $(+)$-path
in $A(R)\cap S_{n,6n}^{\prime \prime }$ from $[0,n]\times \{6n\}$ to a neighbouring vertex of $x$ 
(we use Higuchi's language although our definition is slightly different).
For a fixed $R\in {\cal R}_{n,4n}$,
we denote by $c(R)$ the ``maximal number of cut points in the middle part of $R$ far enough from each other,''
that is, the cardinality of a maximal subset $M(R)\subset R\cap S_{n,6n}^{\prime \prime }$
for which the following properties hold:
\begin{itemize}
\item every $v\in M(R)$ is a cut point of $R$ in $S_{n,6n}$,
\item for all $v_1,v_2\in M(R)$, $d(v_1,v_2)\geq \sqrt{n}$.
\end{itemize}
We shall next compute a lower bound for (a conditional expectation of) $c(R)$ by using the aforementioned
result by Higuchi.

Proposition 5.1 in \cite{Higcoex} concerning the Ising model on the square lattice essentially states that
if both $(+)$-crossings and dual $(-)$-crossings in the long direction of $4n\times n$ rectangles have
probability bounded away from 0, then for an arbitrary fixed horizontal crossing $R$ in the lowest quarter
of an $n$ by $n$ square $S$, irrespective of what the spins of vertices in and below $R$ are, the expected
number of vertices $v$ in $R$ with a $(+)$-path from a neighbour of $v$ to the top of $S$ is arbitrarily
large for all $n$ large enough.
A careful reading of the proof of this proposition shows that the same method works on the triangular
lattice $\mathbb{T}$.
Moreover, we can take the parallelogram $S_{n,6n}$ instead of a square, consider a horizontal crossing
$R\in {\cal R}_{n,4n}$, condition on the spins of vertices in $D(R)$ instead of $L(R)\cup R$,
require that the $(+)$-path from a neighbour of $v\in R$ to the top of $S_{n,6n}$ be in
$A(R)\cap S_{n,6n}^{\prime \prime }$, and the expected number of special vertices
(which are here cut points of $R$ in $S_{n,6n}$) still goes to infinity as $n\to \infty $. In fact,
using Higuchi's notation in \cite{Higcoex}, we see that
since all the cut points considered in the proof are found inside annuli $A_j^{'''}$
which are at distance at least $\frac{5}{2}\cdot 4^j$ from one another (where only
integers $j$ satisfying $\sqrt{n}\leq 2\cdot 4^j$ are considered --- see (5.21) in \cite{Higuchi}), all
cut points considered are automatically at distance at least $\frac{5}{4}\cdot \sqrt{n}$ from one another.
Therefore, if $\mathbb{E}_{\beta}$ denotes the expected value w.r.t.\ $\mu _{\beta}$,
and ${\cal F}_V$ denotes the $\sigma $-algebra generated by $\{\sigma (x):x\in V\}$, we have the following
result.

\begin{Proposition} \label{Hig}
Let $\beta <\beta _c$ and assume that there exists $\delta >0$ such that
\begin{equation}\label{Higcond}
\min \{\mu _{\beta }(H^+ _{3n,n}),\mu _{\beta }(H^{-} _{3n,n})\}\geq \delta
\end{equation}
for every $n\geq 1$. Then we have
$$
\lim _{n\to \infty }\inf _{R\in {\cal R}_{n,4n}}\inf _{E\in {\cal F}_{D(R)}}\mathbb{E}_{\beta }\
(c(R)\mid E)=\infty .
$$
\end{Proposition}

Due to the self-matching property of $\mathbb{T}$ and the $+/-$ symmetry of the model,
for any $n\in \mathbb{N}$, we have
\begin{equation}
\mu _{\beta }(H^{+} _{n,n})= 1/2.
\end{equation}
It follows from this observation and the RSW-type results in \cite{HiguchiRSW}
(which apply to $\mathbb T$ as well as to the square lattice) that condition (\ref{Higcond})
in Proposition \ref{Hig} is satisfied with a proper choice of $\delta $. Furthermore, since
$r=1/2$ corresponds to the Ising model, for all $\beta <\beta _c$, we have
\begin{equation}\label{manycutpoints12}
\lim _{n\to \infty }\inf _{R\in {\cal R}_{n,4n}}\inf _{E\in {\cal F}_{D(R)}}\
\mathbb{E}_{\beta ,1/2}(c(R)\mid E)=\infty .
\end{equation}

\section{Domination lemmas}\label{domination-lemmas}

\subsection{Strategy of the proof of Theorem \ref{rc12}}\label{strategy-of-proof}

In order to motivate the technical results in this section, we give an informal
(and somewhat imprecise) overview of the main step in the proof of Theorem \ref{rc12},
namely the proof that $r_c(\beta) \leq 1/2$. The structure of our proof of this fact is
based on Russo's formulation \cite{Russocrpr} of Kesten's celebrated proof~\cite{Kesten}
of the analogous statement for bond percolation on the square lattice. The proof proceeds
by contradiction, assuming that $r_c(\beta)>1/2$ and showing that this implies the existence
of some $\varepsilon>0$ such that, $\forall r \in [1/2,1/2+\varepsilon]$, the number of FK
clusters which are pivotal for the event corresponding to the presence of a ($-$)-crossing
in a sufficiently large parallelogram $S_{N,6N}$ is very large (in expectation).
By Russo's formula, the expected number of pivotal FK clusters equals the
derivative of the probability of the crossing event.
This leads to a contradiction since the probability of any event has to
remain between 0 and 1, and so its derivative cannot be too large on an interval.

We show in Section \ref{bondmainproofsec} that if we take $\beta <\beta _c$ and assume $r_c(\beta )>1/2$,
then the probability of a horizontal ($-$)-crossing in the lower half $S_{N,3N}$ of the parallelogram
$S_{N,6N}$ is bounded away from 0, uniformly for every $r\in [1/2,r_c(\beta ))$.
We take $r_0$ in that range and consider the lowest such crossing
$R$ and the union $U_R$ of FK clusters of vertices in and below $R$, which is surrounded by a closed
barrier $B$. Since $\beta<\beta_c$, the FK clusters ``tend to be small.'' Therefore, with high probability,
every edge of $B$ is at most at distance $N^{1/4}$ from the set of vertices in and below $R$. Assuming
that this is the case, the vertex boundary of $int(B)$ contains exactly one horizontal crossing of
$S_{N,4N}$, which we call $\Gamma _{B}$. Since the vertices in $\Gamma _{B}\cap S_{N,6N}^{\prime \prime }$
(i.e., the middle part of $\Gamma _{B}$) are in FK clusters of vertices in the lowest horizontal
($-$)-crossing $R$, if $v\in \Gamma _{B}\cap S_{N,6N}^{\prime \prime }$ is a cut point of $\Gamma _{B}$ in
$S_{N,6N}$, then $C^{FK}_v$ is pivotal for $H^- _{N,6N}$. Therefore, from this point on, our goal is
to find a large number of cut points of $\Gamma _{B}$ in $S_{N,6N}$ in $\Gamma _{B}\cap S_{N,6N}^{\prime \prime }$.

In Section \ref{Isingsub}, we used Higuchi's results and the Edwards-Sokal coupling to obtain
equation (\ref{manycutpoints12}), which informally states that for $\beta <\beta _c$ and $r=1/2$,
for any horizontal crossing of a sufficiently large parallelogram, regardless of the values of
the spins of vertices in and below the crossing, the expected number of cut points of the crossing
is arbitrarily large. We would like to use this result to conclude that there are many cut points
of $\Gamma _{B}$ in $S_{N,6N}$ in $\Gamma _{B}\cap S_{N,6N}^{\prime \prime }$. We couple the $r=r_0$ and
the $r=1/2$ case by taking the same random-cluster configuration in $ext(B)$ (which is allowed since
$B$ is a closed barrier), and assigning spins to the FK clusters as follows. We take i.i.d.\ random
variables $(V(C^{FK}_v):v\in ext(B))$ with uniform distribution on the interval $[0,1]$, and assign
spin $+1$ to $v$ if $V(C^{FK}_v)$ is smaller than $r_0$ or $1/2$, respectively. Then, every vertex
which is a cut point in the $r=1/2$ case is a cut point in the $r=r_0$ $(>1/2)$ case as well, since
being a cut point requires the presence of ($+$)-paths only, and every vertex in $ext(B)$ whose spin
is $+1$ at $1/2$ has a $+1$ spin also at $r_0$.

We now would like to use~(\ref{manycutpoints12}), but we cannot do that immediately because at this
point of the proof we have information on the FK clusters of vertices in and below $R$, and not only
on spin values, as required by~(\ref{manycutpoints12}). To circumvent this problem, we will use the
presence of the closed barrier $B$ to show that having information on the FK clusters of vertices in
and below $R$ does not create problems. This is intuitively not surprising, but proving such a result
requires a considerable amount of work, to which the rest of the present section is dedicated.

The proof of Theorem \ref{rc12} can be finished from here as follows. First of all, it follows from
Lemma \ref{conditional-independence} that turning the spin of every vertex in $int(B)$ to $-1$ does
not change the expected number of cut points in $\Gamma _{B}\cap S_{N,6N}^{\prime \prime }$. Then,
Corollary~\ref{barrinwh} implies that this expected number is bounded below by the expected number
without conditioning on $B$ being closed. For the latter expected number, we can use~(\ref{manycutpoints12})
to conclude that the expected number of cut points in $\Gamma _{B}\cap S_{N,6N}^{\prime \prime }$ becomes
arbitrarily large as the size of the parallelogram increases, leading to the desired contradiction,
as discussed earlier.

\subsection{A barrier around $-1$ spins}\label{a-barrier-around-white-vertices}

Our goal in this section is to prove Corollary \ref{barrinwh}. We do this through three lemmas,
using ideas from \cite{BCM} and \cite{KW}. We need a property of the random-cluster measure
$\nu _{p,q}$ on $\mathbb{T}$ from \cite{FKG} (see also \cite{grimmett2}), namely that for all
$q\geq 1$, the so-called ``FKG lattice condition'' holds for $\nu _{p,q}$. We use the following
version of it: for any $E\subset {\cal E}_{\mathbb{T}}$, $e\in {\cal E}_{\mathbb{T}}\setminus E$,
and $\psi ,\zeta \in \{0,1\}^E$ with $\zeta \geq \psi $ (coordinate-wise), we have
\begin{equation}\label{strongFKG}
\nu _{p,q}(\eta (e)=1\mid \eta \equiv \zeta \textrm{ on }E)\geq \nu _{p,q}(\eta (e)=1\mid \eta \equiv \psi \textrm{ on }E).
\end{equation}
This property will play an important role in the following proofs. We state the following lemmas for the measure $\mathbb{P}_{\beta ,r}$
but in fact all statements in this section hold for all DaC measures obtained by replacing $\nu _{p,2}$ in the
construction of $\mathbb{P}_{\beta ,r}$ by $\nu _{p,q}$ with $q\geq 1$.

Inequality (\ref{strongFKG}) informally states that the more edges in a certain set $E$ are open, the more likely
it is that other edges are open as well. The next lemma states that further conditioning on the left hand side on
the event $I$ that the vertices of a certain set $V$ all have the same spin $\kappa $ leaves the inequality
unchanged.

Let $V=\{v_1,v_2,\ldots ,v_k\}\subset {\cal V}_{\mathbb{T}}$ be a set of vertices, $\kappa \in \{-1,+1\}$ a spin value,
$E=\{e_1,e_2,\ldots e_{\ell}\}\subset {\cal E}_{\mathbb{T}}$ a set of edges, $s_1,s_2,\ldots ,s_{\ell}\in \{0,1\}$
and $g_1,g_2,\ldots ,g_{\ell}\in \{0,1\}$ states, with $g_i\geq s_i$ for all $i$. Consider the events
$I=\bigcap _{i=1} ^k \{\sigma (v_i)=\kappa \}$, $A_s=\bigcap _{j=1} ^{\ell} \{\eta (e_j)=s_j \}$,
$A_g=\bigcap _{j=1} ^{\ell} \{\eta (e_j)=g_j \}$
(the case $E=\emptyset ,A_s=A_g=\tilde\Omega $ is also allowed).

\begin{Lemma}\label{edgedom}
For all $e\in {\cal E}_{\mathbb{T}}\setminus E$, we have
\begin{equation}\label{star}
\mathbb{P}_{\beta ,r}(\eta (e)=1\mid A_g,I)\geq \mathbb{P}_{\beta ,r}(\eta (e)=1\mid A_s).
\end{equation}
\end{Lemma}

\noindent {\bf Proof:}
Since
\begin{eqnarray*}
\mathbb{P}_{\beta ,r}(\eta (e)=1\mid A_g,I)& = & \frac{\mathbb{P}_{\beta ,r}(\eta (e)=1,A_g,I)}{\mathbb{P}_{\beta ,r}(A_g,I)}\\
 & = & \frac{\mathbb{P}_{\beta ,r}(I\mid \eta (e)=1,A_g)}{\mathbb{P}_{\beta ,r}(I\mid A_g)}\cdot \mathbb{P}_{\beta ,r}(\eta (e)=1\mid A_g),
\end{eqnarray*}
and $\mathbb{P}_{\beta ,r}(\eta (e)=1\mid A_g)\geq \mathbb{P}_{\beta ,r}(\eta (e)=1\mid A_s)$ by (\ref{strongFKG}),
we have that (\ref{star}) follows from
\begin{equation}\label{star3}
\mathbb{P}_{\beta ,r}(I\mid \eta (e)=1,A_g)\geq \mathbb{P}_{\beta ,r}(I\mid A_g).
\end{equation}
Since
\begin{eqnarray*}
\mathbb{P}_{\beta ,r}(I\mid A_g)& = & \mathbb{P}_{\beta ,r}(I\mid \eta (e)=1,A_g)\mathbb{P}_{\beta ,r}(\eta (e)=1\mid A_g) \\
& + & \mathbb{P}_{\beta ,r}(I\mid \eta (e)=0,A_g)\mathbb{P}_{\beta ,r}(\eta (e)=0\mid A_g),
\end{eqnarray*}
we see that (\ref{star3}) is equivalent to
\begin{equation}\label{star4}
\mathbb{P}_{\beta ,r}(I\mid \eta (e)=1,A_g)\geq \mathbb{P}_{\beta ,r}(I\mid \eta (e)=0, A_g).
\end{equation}

In order to show (\ref{star4}), we will first construct two coupled bond configurations, $\psi _0$
and $\psi _1$, such that $\psi _0$ has distribution $\nu _{p,2}(\cdot \mid \eta (e)=0,A_g)$,
$\psi _1 $ has distribution $\nu _{p,2}(\cdot \mid \eta (e)=1,A_g)$ (both with $p=1-e^{-\beta }$),
and $\psi _0\leq \psi _1$.
Such a coupling can be obtained by setting $\psi _0(e_i)=\psi _1(e_i)=g_i$, $\psi _0(e)=0$,
$\psi _1(e)=1$, then determining the states of the remaining edges one edge at a time in some
deterministic order, using (\ref{strongFKG}) at each step (for a precise way of doing this,
see e.g.\ the proof of Lemma 2 in \cite{KW}).

We could easily finish the proof from here by completing the coupling to obtain
configurations $(\psi _0,\sigma _0)$ and $(\psi _1,\sigma _1)$
with distributions $\mathbb{P}_{\beta ,r}(\cdot \mid \eta (e)=0,A_g)$ and
$\mathbb{P}_{\beta ,r}(\cdot \mid \eta (e)=1,A_g)$
respectively, in such a way that if $I$ occurs in $\sigma _0$, it occurs in $\sigma _1$ as well.
Alternatively, we may notice that given a bond configuration $\psi $, defining $n(\psi )$ as the
number of FK clusters in $\psi $ which contain vertices of $V$, the probability of $I$ is simply
$c^{n(\psi )}$, where $c=r$ if $\kappa =+1$ and $c=1-r$ if $\kappa =-1$. Since $n(\psi _0)\geq n(\psi _1)$
and $0\leq c\leq 1$, this observation concludes proof of (\ref{star4}) and thereby the proof of
Lemma \ref{edgedom}.
\hfill $\Box $\\

Now take $E=\{e_1,e_2,\ldots ,e_{\ell}\},s_1,s_2,\ldots ,s_{\ell},A_s$ as before,
and let $F=\{f_1,f_2,\ldots ,f_m\}\subset {\cal E}_{\mathbb{T}}$ be a set of edges such that
$F\cap E=\emptyset $, and define the event $C(F)=\bigcap _{i=1}^m\{\eta (f_i)=0\}$. Then, as
an easy consequence of (\ref{strongFKG}), we have that for all
$q\geq 1, e\in {\cal E}_{\mathbb{T}}\setminus (E\cup F)$,
$$
\nu _{p,q}(\eta (e)=1\mid A_s)\geq \nu _{p,q}(\eta (e)=1\mid A_s,C(F)).
$$
The next lemma follows from this observation and Lemma \ref{edgedom}.
\begin{Lemma}\label{edgedom2}
For all $e\in {\cal E}_{\mathbb{T}}\setminus (E\cup F)$, we have
$$
\mathbb{P}_{\beta ,r}(\eta (e)=1\mid A_g,I)\geq \mathbb{P}_{\beta ,r}(\eta (e)=1\mid A_s,C(F)).
$$
\end{Lemma}
Note that this statement is still an intuitively clear consequence of (\ref{strongFKG}),
since the additional conditioning on $I$ (i.e.\ that certain vertices all have spin $\kappa $)
on the left hand side of (\ref{strongFKG}) should intuitively increase the probability that
other edges are open, whereas the additional conditioning on $C(F)$ (i.e.\ having even more
edges closed) on the other side should intuitively decrease this probability.

We are now ready to state the main result in this section, which immediately implies the desired Corollary \ref{barrinwh}.
\begin{Lemma}\label{barrin}
Let $V=\{v_1,v_2,\ldots ,v_k\}\subset {\cal V}_{\mathbb{T}}$ be a connected set of vertices,
and take its edge boundary $B=\Delta V=\{f_1,f_2,\ldots ,f_m\}\subset {\cal E}_{\mathbb{T}}$
(which is a barrier).
Consider the events $I=\bigcap _{i=1}^k\{\sigma (v_i)=-1\}$, $C(B)=\bigcap _{j=1}^m\{\eta (f_j)=0\}$,
and let $D$ be an increasing event. Then we have
\begin{equation}\label{barrinequation}
\mathbb{P}_{\beta ,r}(D\mid C(B))\geq \mathbb{P}_{\beta ,r}(D\mid I).
\end{equation}
\end{Lemma}

\noindent {\bf Proof.}
We prove (\ref{barrinequation}) by constructing two coupled realisations
$(\psi _{C(B)},\sigma _{C(B)})$ and $(\psi _I,\sigma _I)$ with distributions
$\mathbb{P}_{\beta ,r}(\cdot \mid C(B))$ and $\mathbb{P}_{\beta ,r}(\cdot \mid I)$
respectively, in such a way that if $D$ occurs in $\sigma _I$, it occurs in
$\sigma _{C(B)}$ as well.

First, we construct the bond configurations $\psi _{C(B)}$ and $\psi _I$ one edge
at a time, using Lemma \ref{edgedom2} at each step, as follows. Fix a deterministic
order of edges in ${\cal E}_{\mathbb{T}}$ starting with edges incident on $v_1,v_2,\ldots ,v_k$.
Take a collection $(U(e):e\in {\cal E}_{\mathbb{T}})$ of i.i.d.\ random variables
having uniform distribution on the interval $[0,1]$. We start with a situation where
$\psi _{C(B)}(e)$ and $\psi _I(e)$ are undetermined for every edge, and determine the
states of edges by the following iteration. We take the first edge in the deterministic
order, and denote it by $e_1$. We declare $\psi _{C(B)}(e_1)=1$ if and only if
$U(e_1)\leq \mathbb{P} _{\beta ,r}(\eta (e_1)=1 \mid C(B))$, and $\psi _I (e_1)=1$ if
and only if $U(e_1)\leq \mathbb{P} _{\beta ,r}(\eta (e_1)=1 \mid I)$. Note that by
Lemma \ref{edgedom2}, $\psi _{C(B)}(e_1)\leq \psi _I(e_1)$.

Let us now assume that the states of $e_1,e_2,\ldots ,e_j$ are determined and
$\psi _{C(B)}(e_i)\leq \psi _I(e_i)$ for $i=1,2,\ldots ,j$. The next edge
$e_{j+1}$ is the next undetermined edge in our deterministic order that shares
a vertex with an edge which is open in $\psi _I$. If no such edge exists,
we simply take the next undetermined edge.

Having chosen $e_{j+1}$, we determine its state by defining
$\psi _{C(B)}(e_{j+1})=1$ if and only if
$U(e_{j+1})\leq \mathbb{P} _{\beta ,r}(\eta (e_{j+1})=1 \mid C(B),\cap _{i=1}^{j}
\{ \eta (e_i)=\psi _{C(B)}(e_i) \} )$
(otherwise we assign $\psi _{C(B)}(e_{j+1})=0$),
and $\psi _I (e_{j+1})=1$ if and only if
$U(e_{j+1})\leq \mathbb{P} _{\beta ,r}(\eta (e_{j+1})=1 \mid I,\cap_{i=1}^{j} \{ \eta (e_i)=\psi _{I}(e_i) \} )$
(otherwise $\psi _I (e_{j+1})=0$).
By the hypothesis $\psi _{C(B)}(e_i)\leq \psi _I(e_i)$ for $i=1,2,\ldots ,j$ and Lemma \ref{edgedom2}, we have that
$\psi _{C(B)}(e_{j+1})\leq \psi _I(e_{j+1})$.

In this way, we obtain bond configurations $\psi _{C(B)}$ with distribution $\mathbb{P} _{\beta ,r}(\cdot \mid C(B))$
and $\psi _I$ with distribution $\mathbb{P} _{\beta ,r}(\cdot \mid I)$ such that $\psi _{C(B)}\leq \psi _I$.
Let us fix $j^*$ to be the index of the last edge chosen by the iteration which is connected by a $\psi _I$-open
edge path to any of the vertices $v_1,v_2,\ldots ,v_k$.
The first part of the iteration (i.e.\ before $e_{j^*+1}$ is chosen)
``explores'' the FK clusters in $\psi _I$ of the vertices
$v_1,v_2,\ldots ,v_k$, and when it ends, $V$ is surrounded by a barrier $B_2$ (which consists of edges
from $e_1,e_2,\ldots ,e_{j^*}$) which is closed in $\psi _I$.
Since $\psi _I \geq \psi _{C(B)}$, $B_2$ is closed in $\psi _{C(B)}$ as well.
Using Lemma \ref{conditional-independence}, we obtain
\begin{eqnarray*}
& \mathbb{P} _{\beta ,r}(\eta (e_{j^*+1})=1 \mid C(B),\cap_{i=1}^{j^*} \{ \eta(e_i)=\psi_{C(B)}(e_i) \} ) &\\
& =\mathbb{P} _{\beta ,r}(\eta (e_{j^*+1})=1 \mid I,\cap_{i=1}^{j^*} \{ \eta(e_i)=\psi_{I}(e_i) \} ),&
\end{eqnarray*}
which implies $\psi _{C(B)}(e_{j^*+1}) =\psi _I(e_{j^*+1})$.
Using the same argument, it is easy to prove by induction that the remaining part of the iteration yields $\psi _{C(B)} =\psi _I$ in $ext(B_2)$.

We now define the spin configuration $\sigma _I$ by assigning $+1$ with probability $r$, $-1$ with
probability $1-r$ independently to the $\psi _I$ FK clusters in $ext(B_2)$ (according to some
deterministic order), and assigning $\sigma _I(v)=-1$ to each $v\in int(B_2)$.
This gives the correct distribution since every vertex in $int(B_2)$ is in the same FK cluster as
one of the vertices $v_1,v_2,\ldots ,v_k$.
We finish the coupling by defining $\sigma _{C(B)}$ in the following way.
We assign $+1$ with probability $r$, $-1$ with probability $1-r$ independently to the $\psi_{C(B)}$
FK clusters in $int(B_2)$ (according to some deterministic order), and define
$\sigma_{C(B)}(v)=\sigma _I(v)$ for all $v\in ext(B_2)$ (since $\psi _{C(B)} =\psi _I$ in $ext(B_2)$,
we get the right distribution).
Let us assume that $D$ occurs in $\sigma _I$. It is important to notice that all vertices that
have spin $+1$ in $\sigma _I$ are in $ext(B_2)$, where $\sigma _{C(B)}=\sigma _I$, so they have
spin $+1$ also in $\sigma _{C(B)}$.
Since $D$ is an increasing event, this observation shows that $D$ occurs in $\sigma _{C(B)}$ as well.
This concludes the proof of Lemma \ref{barrin}.
\hfill $\Box $

\begin{Corollary}\label{barrinwh}
If $V=\{v_1,v_2,\ldots ,v_k\}\subset {\cal V}_{\mathbb{T}}$ is a connected set of vertices,
$B=\Delta V=\{f_1,f_2,\ldots ,f_m\}\subset {\cal E}_{\mathbb{T}}$ is its edge boundary,
and we consider the events $I=\bigcap _{i=1}^k\{\sigma (v_i)=-1\}$, $C(B)=\bigcap _{j=1}^m\{Y(f_j)=0\}$,
and an increasing event $D\in {\cal F}_{ext(B)}$, then we have
\begin{equation}\label{barrinwheq}
\mathbb{P}_{\beta ,r}(D\mid C(B),I)\geq \mathbb{P}_{\beta ,r}(D\mid I).
\end{equation}
\end{Corollary}

\noindent {\bf Proof.}
Since $B$ is a barrier, $I\in {\cal F}_{int(B)}$, and $D\in {\cal F}_{ext(B)}$, we have by Lemma \ref{conditional-independence} that
$\mathbb{P}_{\beta ,r}(D\mid C(B),I)=\mathbb{P}_{\beta ,r}(D\mid C(B))$.
Therefore, Lemma \ref{barrin} gives the statement.
\hfill $\Box $

\section{Proofs of Theorems \ref{rc12}--\ref{contbond}}\label{bondmainproofsec}

In order to prove $r_c(\beta )=1/2$ in Theorem \ref{rc12}, we only need
to show $r_c(\beta )\leq 1/2$, since $r_c(\beta )\geq 1/2$ is implied by
Proposition 1.8 in \cite{BCM}. By Theorem \ref{implperc}, it suffices to
prove that $\limsup_{n \to \infty} {\mathbb P}_r(V^+_{n,3n}) = 1$ when
$r=1/2+\varepsilon$ for all $\varepsilon>0$. We shall prove that
the assumption of the contrary implies the presence of too many pivotal FK
clusters for a certain event, leading to a contradiction.
(For a more detailed summary of the proof, see Section \ref{strategy-of-proof}.)

\begin{Theorem}\label{q2rc}
For any $\beta <\beta _c$ and $\varepsilon >0$, we have that
$$
\limsup \limits _{n\to \infty }\mathbb{P}_{\beta ,1/2+\varepsilon }(V^+ _{n,3n})=1
$$
\end{Theorem}

\noindent \textbf{Proof: }Let us assume that there exist $\beta <\beta _c,\varepsilon >0$ such that
\begin{equation}\label{assumption}
\limsup \limits _{n\to \infty }\mathbb{P}_{\beta ,1/2+\varepsilon }(V^+ _{n,3n})<1,
\end{equation}
and fix such a $\beta $ and $\varepsilon $. We shall derive a contradiction from (\ref{assumption}).
Due to the self-matching property of $\mathbb{T}$, (\ref{assumption})
implies that there exists $\gamma >0$ such that for all $n$ large enough,
\begin{equation}\label{Nlarge1}
\mathbb{P}_{\beta ,1/2+\varepsilon }(H^- _{n,3n})>\gamma .
\end{equation}
By (\ref{Nlarge1}), monotonicity, (\ref{manycutpoints12}), and elementary properties of the
exponential function, it is possible to choose an integer $N$ large enough so that for $n\geq N$,
the following inequalities hold:
\begin{eqnarray}
\mathbb{P}_{\beta ,r}(H^- _{n,3n})& > & \gamma \,\,\,\,\,\,\,\,\, \forall r \in [1/2,1/2+\varepsilon ],
\label{horcross}\\
\inf _{R\in {\cal R}_{n,4n}}\inf _{E\in {\cal F}_{D(R)}}\mathbb{E}_{\beta ,1/2}(c(R)\mid E) & > &
\frac{2}{\varepsilon \gamma }, \label{manycut}\\
(n+1)(6n+1)e^{-n^{1/4}\psi (p)} & < & \frac{\gamma }{2},\label{shortranges}
\end{eqnarray}
where $\psi (p)$ is the same as in Theorem \ref{lexpdecay2}.
Fix such an $N$ and an arbitrary $r_0 \in [1/2,1/2+\varepsilon ]$. We shall show that,
denoting the number of FK clusters which are pivotal for $H^-_{N,6N}$ by $n(H^-_{N,6N})$,
we have
\begin{equation}\label{manypiv2}
\mathbb{E}_{\beta ,r_0}(n(H^- _{N,6N}))>\frac{1}{\varepsilon }.
\end{equation}

For $R\in {\cal R}_{N,3N}$, we define
\begin{eqnarray*}
B(R) & = & \{ B \subset {\cal E}_{\mathbb T}: B \textrm{ is a barrier}; \partial (L(R)\cup R)\subset int(B);\\
 &  & \forall \hspace{0.05cm}e\in B \hspace{0.1cm}, d(e,\partial (L(R)\cup R))\leq N^{1/4};\\
 &  & \partial \hspace{0.02cm}int(B) \textrm{ contains exactly one}\\
 &  & \textrm{horizontal crossing of }S_{N,4N}\}.
\end{eqnarray*}
(The motivation for this definition is that since $\beta <\beta _c$, FK clusters are small,
hence with high probability, the ``tightest'' closed barrier surrounding $L(R)\cup R$ is
contained in $B(R)$.)
For $B\in B(R)$, we denote the horizontal crossing of $S_{N,4N}$ contained in
$\partial \hspace{0.02cm} int(B)$ by $\Gamma _{B}$.
We also define $l(R)$ to be the event that $R$ is the lowest horizontal ($-$)-crossing in $S_{N,6N}$.
For $R\in {\cal R}_{N,3N}$, $B\in B(R)$, we denote the union of FK clusters
$\bigcup _{v\in L(R)\cup R}C^{FK}_v$ by $U_R$, the event 
$\bigcap _{v\in L(R)\cup R}\{\mathcal{D}(v)\leq N^{1/4}\}$ by $t(R)$,
and consider the event
$$Q(R,B)=\{l(R)\}\cap \{B=\Delta U_R\}\cap \{t(R)\}.$$
Then we obtain
$$
\mathbb{E}_{\beta ,r_0}(n(H^- _{N,6N})) \geq \sum _{R\in {\cal R}_{N,3N}}\
\sum _{B\in B(R)}\mathbb{E}_{\beta ,r_0}(n(H^- _{N,6N})\mid Q(R,B))\mathbb{P}_{\beta ,r_0}(Q(R,B))
$$
\begin{equation}\label{firstest}
\geq \sum _{R\in {\cal R}_{N,3N}}\sum _{B\in B(R)}\mathbb{E}_{\beta ,r_0}(c(\Gamma _{B})\mid Q(R,B))
\mathbb{P}_{\beta ,r_0}(Q(R,B)),
\end{equation}
where the second inequality follows from a pointwise comparison: conditioned on $Q(R,B)$, we have
$n(H^- _{N,6N})\geq c(\Gamma _{B})$, due to the following reasons. Using the notation from the
definition of $c(\Gamma _{B})$ (see Section \ref{Isingsub}), conditioned on $Q(R,B)$, the FK cluster
of every vertex $v$ in $M(\Gamma _{B})$ is pivotal for $H^- _{N,6N}$ since $v$ is a cut point of
$\Gamma_{B}$ in $S_{N,6N}$, and $R$ is the lowest horizontal ($-$)-crossing in $S_{N,6N}$. It is
important to note that every $v\in M(\Gamma _{B})$ is indeed in the FK cluster of a vertex in $R$
(i.e., of a vertex in the lowest horizontal ($-$)-crossing), not of a vertex in $L(R)$ (there is
no other possibility due to $\{B=\Delta U_R\}$). This is the case since
$M(\Gamma _{B})\subset \Gamma _{B}\cap S_{N,6N}^{\prime \prime }$ --- since none of the vertices below
$R$ has a dependence range larger than $N^{1/4}$, none of the FK clusters of the vertices in $L(R)$
is large enough to go around $R$ and reach the middle part $S_{N,6N}^{\prime \prime }$ of the parallelogram
$S_{N,6N}$. The last step necessary for proving the conditional pointwise comparison is to notice
that for $v_1,v_2\in  M(\Gamma _{B}), v_1\neq v_2$, we have $C_{v_1}^{FK}\neq C_{v_2}^{FK}$ since
$d(v_1,v_2)\geq \sqrt{N}$ and, conditioned on $Q(R,B)$, none of the vertices in $L(R)\cup R$ has
a dependence range greater than $N^{1/4}$. Therefore, different vertices in $M(\Gamma _{B})$ belong
to different pivotal FK clusters.

The next step is to give a lower bound for the expectation via a comparison with the case with parameter
$r=1/2$. We shall first work with probabilities, then we will sum them up to get back the expectation.
Let us denote $(N+1)(6N+1)$ (i.e.\ the number of vertices in $S_{N,6N}$) by $K$.
For a barrier $B$, we define the events $C(B)=\bigcap _{e\in B} \{\eta (e)=0\}$ and
$W(B)=\bigcap _{v\in int(B)} \{\sigma (v)=-1\}$.
Since for every $R\in {\cal R}_{N,3N}$, $B\in B(R)$, $i=1,\ldots ,K$,
we have $\{c(\Gamma _{B})\geq i\}\in \mathcal{F}_{A(\Gamma _B)\cap
S_{N,6N}^{\prime \prime }}\subset \mathcal{F}_{ext(B)}$,
$\{l(R)=R\}\in \mathcal{F}_{L(R)\cup R}\subset \mathcal{F}_{int(B)}$,
$W(B)\in \mathcal{F}_{int(B)}$, and the event
$\{B=\Delta U_R\}\cap \{t(R)\}$
depends on the state of edges in $int(B)\cup B$ only, it follows from a repeated use of
Lemma \ref{conditional-independence} that for all $R,B$, and $i$, we have
\begin{eqnarray}
\mathbb{P}_{\beta ,r_0}(c(\Gamma _{B})\geq i\mid Q(R,B)) & = & \mathbb{P}_{\beta ,r_0}(c(\Gamma _{B})
\geq i\mid C(B))\nonumber \\
 & = & \mathbb{P}_{\beta ,r_0}(c(\Gamma _{B})\geq i\mid C(B),W(B)).\label{qrbtocbwb}
\end{eqnarray}
Coupling the measures with $r=r_0$ and $r=1/2$ by taking the same bond configurations in $ext(B)$
(see Section \ref{strategy-of-proof}), we see that
\begin{equation}\label{rto12}
\mathbb{P}_{\beta ,r_0}(c(\Gamma _{B})\geq i\mid C(B),W(B)) \geq
\mathbb{P}_{\beta ,1/2}(c(\Gamma _{B})\geq i\mid C(B),W(B)).
\end{equation}
Since for all $i$, $\{c(\Gamma _{B})\geq i\}\in \mathcal{F}_{ext(B)}$
is an increasing event, we can use Corollary \ref{barrinwh} to conclude that
\begin{equation}\label{cbwbtowb}
\mathbb{P}_{\beta ,1/2}(c(\Gamma _{B})\geq i\mid C(B),W(B))\geq \mathbb{P}_{\beta ,1/2}(c(\Gamma _{B})\geq i\mid W(B)).
\end{equation}
Summing up for $i=1,\ldots ,K$, using (\ref{qrbtocbwb}),(\ref{rto12}),(\ref{cbwbtowb}) and then (\ref{manycut}), we obtain that
for every $R\in {\cal R}_{N,3N}$, $B\in B(R)$, a.s.,
\begin{eqnarray}
\mathbb{E}_{\beta ,r_0}(c(\Gamma _{B})\mid Q(R,B)) & = & \sum _{i=1} ^{K}\mathbb{P}_{\beta ,r_0}(c(\Gamma _{B})\geq i\mid Q(R,B))\nonumber \\
 & \geq & \sum _{i=1} ^{K}\mathbb{P}_{\beta ,1/2}(c(\Gamma _{B})\geq i\mid W(B))\nonumber \\
 & = & \mathbb{E}_{\beta ,1/2}(c(\Gamma _{B})\mid W(B))\nonumber \\
 & > & \frac{2}{\varepsilon \gamma }. \label{eq4}
\end{eqnarray}

Finally we need to note that for a crossing $R\in {\cal R}_{N,3N}$, if 
$t(R)$ occurs, then
$\Delta U_R\in B(R)$. Therefore,
\begin{eqnarray*}
\sum _{R\in {\cal R}_{N,3N}}\sum _{B\in B(R)}\mathbb{P}_{\beta ,r_0}(Q(R,B)) & = & \sum _{R\in {\cal R}_{N,3N}} \
\mathbb{P}_{\beta ,r_0}(l(R)\cap t(R))\\
 & \geq & \mathbb{P}_{\beta ,r_0}(H^- _{N,3N})-\mathbb{P}_{\beta ,r_0}(\bigcup _{v\in R_{n,6n}}\mathcal{D}(v)> N^{1/4})\\
 & \geq & \gamma -(N+1)(6N+1) \nu _{p,2}(\mathcal{D}(0)> N^{1/4})\\
 & \geq & \gamma -(N+1)(6N+1)e^{-N^{1/4}\psi (p)}\\
 & \geq & \gamma /2,
\end{eqnarray*}
where we used the translation invariance of $\nu _{p,2}$, (\ref{horcross}), Theorem \ref{lexpdecay2},
and (\ref{shortranges}).
Using (\ref{firstest}), (\ref{eq4}), and this computation, we obtain that
\begin{eqnarray*}
\mathbb{E}_{\beta ,r_0}(n(H^- _{N,6N})) & > & \sum _{R\in {\cal R}_{N,3N}}\
\sum _{B\in B(R)}\frac{2}{\varepsilon \gamma }\hspace{0.1cm} \mathbb{P}_{\beta ,r_0}(Q(R,B))\\
 & \geq & \frac{2}{\varepsilon \gamma }\cdot \frac{\gamma }{2}=\frac{1}{\varepsilon },
\end{eqnarray*}
as desired.

Since (\ref{manypiv2}) can be proved for all $r\in [1/2,1/2+\varepsilon ]$
with the same method, we obtain by Theorem \ref{Russocons} that
\begin{displaymath}
\sup \limits _{r\in \left[ 1/2,1/2+\varepsilon \right]}\frac{d}{dr}
\mathbb{P}_{\beta ,r}(H_{N,3N}^{-})<-\frac{1}{\varepsilon },
\end{displaymath}
which leads to a contradiction since it yields
\begin{eqnarray}
\mathbb{P}_{\beta ,1/2+\varepsilon }\left( H_{N,3N}^{-}\right) & \leq &
\mathbb{P}_{\beta ,1/2}\left( H_{N,3N}^{-}\right) \
+\varepsilon \sup \limits _{r\in \left[ 1/2,1/2+\varepsilon \right]}\frac{d}{dr}
\mathbb{P}_{\beta ,r}(H_{N,3N}^{-}) \nonumber \\
& < & \mathbb{P}_{\beta ,1/2}\left( H_{N,3N}^{-}\right) - 1. \nonumber
\end{eqnarray}
\hfill $\Box $\\

\noindent \textbf{Sketch of the Proof of Theorems \ref{rc12} and \ref{contbond}.}
As remarked at the beginning of this section, for all $\beta <\beta _c$,
$r_c(\beta )\geq 1/2$ follows from Proposition 1.8 of \cite{BCM}, and
$r_c(\beta )\leq 1/2$ from Theorem \ref{q2rc} and Theorem \ref{implperc}. Hence,
$r_c(\beta )=1/2$. The exponential tail of the distribution of the size of the
($+$)-cluster of the origin for $r<1/2$ can be proved similarly to Theorem 2 in \cite{Voronoi}.
The statement concerning the critical case $r=1/2$ has been proved in Proposition 1.8
of \cite{BCM}. For $\beta <\beta_c$, the ergodicity of $\mathbb{P}_{\beta ,r}$
(which follows from the ergodicity of $\nu _{p,2}$) guarantees the presence of an
infinite ($+$)-cluster when $r>1/2$.
The uniqueness of the infinite ($+$)-cluster follows from a result in \cite{BuKe},
which implies that if a probability measure
$\mu $ on $\{-1,+1\}^{\cal V_{\mathbb{T}}}$ is translation invariant and satisfies
the finite energy condition \cite{NSch}, then $\mu $-a.s.\ there exists at most one
infinite cluster of $+1$'s. If $\beta <\infty$ and $0<r<1$, then the spin marginal
of $\mathbb{P}_{\beta ,r}$ clearly satisfies both properties.

Theorem \ref{contbond} about the continuity of $\Theta (\beta ,r)$ in $r$ for $\beta <\beta _c$
follows from $\Theta (\beta ,1/2)=0$ and the uniqueness of the infinite ($+$)-cluster by standard
methods (see~\cite{vdbk}), in the same way as the analogous result in \cite{BCM}.
\hfill $\Box $\\

\begin{Remark}
\emph{In all the proofs in this paper, the FKG inequality and RSW-type arguments
are used for $\mathbb{P}_{\beta ,r}$ only \emph{at} the critical point $r=1/2$,
never away from it. This way of proving classical percolation results can be useful
in the case of models, like the present one, where the (conjectured) critical point has
special properties and is better understood compared to other values of the parameter.}
\end{Remark}

\medskip\noindent
{\bf Acknowledgements} We thank Rongfeng Sun for drawing our attention to \cite{KW}.
F.C. thanks Reda J\"urg Messikh and Akira Sakai for interesting discussions at an early
stage of this work.

\end{document}